\documentclass[12pt]{article}

\usepackage[english]{babel}
\usepackage{amssymb}
\usepackage{latexsym}
\usepackage[round]{natbib}
\usepackage{url}
\usepackage{color}

\renewcommand{\Re}{\mathbb{R}}

\newcommand{\prob}{\mathrm{Pr}}
\newcommand{\vat}{\mathbb{E}}

\newcommand\cvd{\hfill\raisebox{-0.5ex}{$\Box$}}

\newcommand{\bc}{\begin{center}}
\newcommand{\ec}{\end{center}}
\newcommand{\bitem}{\begin{itemize}}
\newcommand{\eitem}{\end{itemize}}
\newcommand{\be}{\begin{eqnarray*}}
\newcommand{\ee}{\end{eqnarray*}}
\newcommand{\ben}{\begin{eqnarray}}
\newcommand{\een}{\end{eqnarray}}
\usepackage{color}

\newcommand\by{\mathbf{y}}

\newenvironment{proof}[1][Proof]{\noindent\textbf{#1.} }{\ \rule{0.5em}{0.5em}}
\newcommand\limn{\lim_{n \rightarrow \infty}}
\newcommand\Pntheta{P_{n,\theta}}
\newcommand\Pnthetan{P_{n,\theta_n}}
\newcommand\Pnthetao{P_{n,\theta_0}}
\newcommand\Pnbeta{P_{n,\beta}}

\newcommand\Pnbetan{P_{n,\beta_n}}

\newcommand\Pno{P_{n,0}}
\newcommand\Pnrsqrtn{P_{n,r/\sqrt{n}}}

\newcommand\pintheta{\pi_{n}(\theta)}
\newcommand\pinthetao{\pi_{n}(\theta_0)}

\newcommand\pinbetan{\pi_{n}(\beta_n)}

\newcommand\sommaxisq{\sum_{i=1}^{n}x_i^2}
\newcommand\sqrtsommaxisq{\sqrt{\sum_{i=1}^{n}x_i^2}}
\newcommand\sommaxiyi{\sum_{i=1}^n x_i y_i}

\usepackage{layout}
\newtheorem{definition}{Definition}[section]

\newtheorem{lemma}{Lemma}[section]
\newtheorem{prop}{Proposition}[section]

\newtheorem{example}{Example}[section]

\begin{document}
\title{Reconciling Model Selection and Prediction}
\author{George Casella
\\
Department of Statistics, University of Florida,
Gainesville, FL 32611.
\\
\texttt{casella@stat.ufl.edu}
\\
\and
Guido Consonni
\\
Dipartimento di Economia Politica e Metodi Quantitativi
\\University of Pavia, 27100 Pavia, Italy.
\\
  \texttt{guido.consonni@unipv.it}
}
\date{\today }
\maketitle

\begin{abstract}
\noindent
It is known that there is a dichotomy in the performance of model selectors.  Those that are consistent (having the ``oracle property") do not achieve the asymptotic minimax rate for prediction error.  We look at this phenomenon closely, and argue that the set of parameters on which this dichotomy occurs is extreme, even pathological, and should not be considered when evaluating model selectors.  We characterize this set, and show that, when such
parameters are dismissed from consideration, consistency and asymptotic minimaxity can be attained simultaneously.
\end{abstract}




\newpage


\section{Introduction}\label{sec:intro}

Model selection is an important area of statistical practice and research. However,  model selection often represents a first step towards our main goal, which may be estimation or prediction. A stylized scheme is the following: first use a model-selection procedure to select a model, and then proceed with inference conditionally on the chosen model.   This method leads to so-called  \lq \lq post-model-selection" estimators
(or predictors). Our concern here is with the asymptotic risk of such procedures.


The review paper \cite{Leeb:Pots:2005} argues that the (data-driven)  model selection step typically has dramatic effect on the sampling properties of the estimators; see also \cite{Leeb:Pots:2006}.  These properties are quite different from their single-model counterpart, and cannot be ignored even when the sample size is large and when the model selector is consistent.

Although interest in the performance of post-model selection estimators has gained momentum over the last years, the problem has been around for a few decades and can be traced backed, in its essence,   to Hodges' estimator.
For an interesting  discussion of Hodges' estimator see
\citet[Example 8.1]{Vaar:1998}.

Within the context of regression functions, and for squared-error loss, \cite{Yang:2005} has shown that consistent model selection procedures, such as BIC, produce estimators which  cannot attain the asymptotic minimax rate.  Failure to attain this optimal rate extends also to model combination, or Bayesian model averaging with subjectively specified 
priors. On the other hand, it is known that AIC, which is inconsistent, does attain the minimax rate; see  \citet [Proposition 1]{Yang:2005}. This tension between model consistency and inference-optimality is sometimes referred to as the \lq \lq AIC-BIC dilemma\rq \rq{}.


Attempts at overcoming this dilemma include  adaptive model selection, which, unlike AIC or BIC, employs  a data-driven penalty to  achieve both consistency and optimality;  for a brief account  see \citet[sec. 1.4]{Yang:2005}.   Recently,  van Erven, Gr\"{u}nwald and Rooij
introduced the notion of \lq \lq switch distribution", as an alternative to standard model selection methods, such as the Bayes factor and leave-one-out cross-validation, in an effort to combine the strengths of AIC and BIC; see their 2008 Technical Report
 entitled Catching Up Faster by Switching Sooner: a Prequential Solution to the AIC-BIC Dilemma
(arXiv:0807.1005v1 [math.ST]).

The remainder of the paper is organized as follows. In Section \ref{sec:some} we review the criteria for asymptotic comparison of tests based on power against local alternatives, and show how it relates to the prediction problem in linear regression.
In Section \ref{sec:conflict} we revisit the result of Yang (2005) for the simple linear regression example. In particular we provide an evaluation of the proof he gave for his Theorem 1, and we link his sequence of alternatives to the categorization  given in Lemma \ref{lem:14.30}.  This puts into perspective, and actually explains why the failure to attain the minimax rate occurs; we also show that the minimax rate is achieved for sequences of type 3  in Lemma \ref{lem:14.30}, which are recognized as the only reasonable sequences for the asymptotic comparison of tests.  In Section \ref{sec:contiguity} we comment on  the use of \textit{contiguity} for proving lack of minimax rate for consistent model selectors.  Finally, Section \ref{sec:concluding} offers some concluding remarks.

\section{Asymptotic Comparison of Tests}
\label{sec:some}
In this section we review some results on asymptotic test comparison.  In particular we are interested in the categorization of local alternatives, and how their convergence to the null value interacts with the power of a test.  We then look at the simple example of model selection in linear regression.

\subsection{Categorizing Alternatives}

Consider a sequence of  statistical models $ \{ \Pntheta , \, \theta \in \Theta   \}$ for observations $\by_n:=(y_1,\ldots,y_n)$, $n=1, 2, \ldots$, where we want to test
$H_0:\theta =\theta_0  \, \textit{vs} \, H_1: \theta >\theta_0$.
%
If $\pintheta$ is the power function of a test, for most reasonable tests it holds true that $\lim_{n \rightarrow \infty} \pintheta = 0$ if $\theta=\theta_0$, and $\limn \pintheta = 1$ if $\theta > \theta_0$. This is to be expected because, with arbitrarily many observations, it should be possible to tell the null and the alternative apart with complete accuracy.  This fact means that to compare tests asymptotically we should make the problem \lq \lq harder\rq \rq{}. One way to do this is to consider  a sequence of testing problems
\ben
\label{eq:H0vsH1}
H_0:\theta =\theta_0  \, \textit{vs} \, H_{1n}: \theta =\theta_n,
\een
where $\theta_n > \theta_0$ and $\theta_n \rightarrow \theta_0^+$ in specific ways.


The $L_1$-distance between $P_{n, \theta_0}$ and $P_{n, \theta_n}$
provides a useful characterization of the power function associated with (\ref{eq:H0vsH1}).
Denote the $L_1$-distance between two probability measures $P$ and $Q$ (having density $p$, respectively $q$, with respect to a common measure $\mu$)  by $||P-Q ||$. Then
$$
||P-Q || \stackrel{def}{=}\int |p-q|d\mu =2 \sup_A |P(A)-Q(A) | \stackrel{def}{=} 2 ||P-Q ||_{TV},
$$
where the second equality follows from the well-known relationship between $L_1$ distance and {\it total variation norm} $||P-Q ||_{TV}$.

The following lemma relates the $L_1$-distance to the power function.
\begin{lemma} \textsc{\citep[lemma 14.30]{Vaar:1998}}
\label{lem:14.30}

The power function $\pi_n$ of any test satisfies
\ben
\label{eq:14.30}
\pintheta-\pinthetao \leq  \frac{1}{2}|| \Pntheta-\Pnthetao   ||.
\een
For any $\theta$ and $\theta_0$ there exists a test whose power function attains equality.
\end{lemma}
The implications of Lemma \ref{lem:14.30} are
\begin{enumerate}
\item
If $||  \Pnthetan-\Pnthetao   || \rightarrow 2$, then the sequence $\theta_n$ is converging to $\theta_0$ at a \textit{slow} rate, so that the
two hypotheses are  \textit{strongly separated}. In this case the \textit{difference} $\pintheta-\pinthetao$ tends to 1, which means that we can get all sort of tests; in particular, since equality can be attained, there exist a sequence of tests with power tending to 1 and size tending to 0 (a \textit{perfect}
sequence of tests).

\item
If $||  \Pnthetan-\Pnthetao   || \rightarrow 0$, then the sequence $\theta_n$ is converging to $\theta_0$ at a \textit{fast} rate, so that the
two hypotheses are \textit{weakly separated}. In this case the power of any sequence of tests is asymptotically less than the level (every sequence of tests is {\it worthless}).

\item

If $||  \Pnthetan-\Pnthetao   ||$ is bounded away from 0 and 2, then the sequence $\theta_n$ is converging to $\theta_0$ at a rate such that the
two hypotheses are  \textit{well separated}. In this case, there exists no perfect sequence of tests, but not every test is worthless either.

\end{enumerate}

The consensus in the literature, see for example   \citet [ sec. 13.1]{Lehm:Roma:2005} or    \citet[sec. 14.5]{Vaar:1998}, is that situation 3 is the only reasonable one for the comparison of tests, otherwise the problem is ``asymptotically degenerate".
For iid observations from smooth models, case 3 occurs when $\theta_n$ converges to $\theta_0$
at rate $1/\sqrt{n}$.

Easier calculation often results when using {\it Hellinger distance} rather than the $L_1$-distance.  The Hellinger distance between $P$ and $Q$ is $ H(P,Q)= \left \{  \int (\sqrt{p}-\sqrt{q})^2 d\mu   \right  \}^{1/2}$, and we can rewrite its square as  $H^2(P,Q)=2-2A(P,Q)$, where $A(P,Q)=\int \sqrt{pq} d\mu $ is called the \textit{Hellinger affinity}.  The following inequality holds \citep{Vaar:1998}:
\begin{equation}
\label{eq:inequalityHTotal}
H^2(P,Q) \leq || P-Q || \leq \min \{ 2-A^2(P,Q),2H(P,Q)  \}.
\end{equation}
\subsection{Example:  Simple Linear Regression}
\label{sec:testing}
Here we look in detail at the simple testing problem considered by Yang (2005).
For $i=1, \ldots,n$ let
\ben \label{eq:simpleregressionmodel}
H_0: \, y_i=\epsilon_i \mbox{ and } H_1: \, y_i = \beta x_i+\epsilon_i; \quad \beta>0, \quad
\een
where $\epsilon_i  \stackrel{iid} \sim N(0,1)$.

Recall that two sequences $a_n$ and $b_n$ are said to be of the same order, written $a_n \asymp b_n$, when
 there exist constants $0<r<R< \infty$ and an integer $n_0$ such that, for $n >n_0$,
  $r<|a_n/b_n |<R$.
We start with a simple lemma whose proof is straightforward.
\begin{lemma}
\label{lem:HellingerAffinity}
Let $\Pno$ be the probability measure associated with $H_0$ and $\Pnbeta$ that associated with  $H_1$ in (\ref{eq:simpleregressionmodel}).
Assume that $ \limn \sum_{i=1}^n x_i^2/n $ is a strictly positive constant, so  that $\sum_{i=1}^n x_i^2 \asymp n$.  The Hellinger affinity is
\ben
\label{eq:HellingerAffinityLinearRegression}
A(\Pno,\Pnbeta)=\exp \left \{ -\frac{\beta^2}{8} \sum_{i=1}^n x_i^2  \right \}.
\een
\end{lemma}
As a direct consequence of this lemma, we can characterize sequences $\beta_n$ as follows.

\begin{enumerate}
\item
if $\beta_n= c_n /\sqrt{n}$, with $c_n \rightarrow \infty$, then
 $|| \Pno - \Pnbetan || \rightarrow 2$  and the two models are strongly separated;
\item
if $\beta_n= c_n /\sqrt{n}$, with $c_n \rightarrow 0^+$, then
 $|| \Pno - \Pnbetan || \rightarrow 0$  and the two models are weakly separated;
\item if $\beta_n \asymp(1/\sqrt{n})$, then
 $|| \Pno - \Pnbetan ||$  is bounded away from 0 and 2, so that the two hypotheses are well separated.
\end{enumerate}

\vspace{.3cm}
Following \citet [p. 498]{Lehm:Roma:2005} we conclude that also for this regression problem the problem of testing $\Pno$ versus $\Pnbetan$ is degenerate unless $\beta_n \asymp 1/\sqrt{n}$.    These are therefore  the only meaningful local alternative sequences for evaluating tests.


\subsection{Consistency and Prediction}\label{sec:consistandpred}

The simple testing problem described in  (\ref{eq:simpleregressionmodel}) can be cast as a model selector by defining $A_n=\{(x_i, y_i), i=1, \ldots n: H_1 \mbox{ is selected }\}$, and estimating $\beta$ with the post-model-selection estimator
$\hat \beta I(A_n)$, where $I(A_n)$ is the indicator function of the set $A_n$.  For simplicity we take $\hat \beta$ to be the least squares estimator under $H_1$.
%

Given $(x_i, y_i), i=1, \ldots n$, we run our model selection procedure and then we predict at values $x_i^\ast, i=1, \ldots,m$, where the $x_i^\ast$ may be the same as the original $x_i$, or not.  The average prediction error is
$$
\frac{1}{m}\sum_{i=1}^m (\beta x_i^\ast-\hat \beta I(A_n)x_i^\ast)^2 =\left(\frac{1}{m}\sum_{i=1}^m  x_i^{\ast \; 2} \right)(\beta -\hat \beta I(A_n))^2,
$$
which shows why it doesn't matter whether we predict $x_i^\ast$ or $x_i$, or how many $x_i^\ast$ we predict (as long as $m/n$ is finite as $n \rightarrow \infty$).
Taking expectations gives the predictive
risk function
\begin{eqnarray*}
R_m (\beta, A_n)&=&\left(\frac{1}{m}\sum_{i=1}^m (x_i^*)^2\right) {\vat}_{\beta} (\beta -\hat \beta I(A_n))^2 \\
&=&\left(\frac{1}{m}\sum_{i=1}^m (x_i^*)^2\right) [{\vat}_{\beta} \{(\beta  -\hat \beta )^2 I(A_n)\} + \beta^2 P_{\beta}(A_n^c)].
\end{eqnarray*}
The predictor attains the asymptotic minimax rate if $ n \sup_\beta R_m (\beta, A_n) \rightarrow$ constant as $n \rightarrow \infty$.
The quantity $nR_m (\beta, A_n)$ is called  the \textit{scaled} risk (function).
Recalling that  $\sum_{i=1}^m (x_i^*)^2 \asymp m $,
 we only need to be concerned with
\begin{eqnarray}
\label{eq:n sup_beta Rm (beta, An)}
n \sup_\beta R_m (\beta, A_n)
\asymp
n [{\vat}_{\beta} \{(\beta  -\hat \beta )^2 I(A_n)\} + \beta^2 P_{\beta}(A_n^c)].
\end{eqnarray}
The first term is bounded by $n {\vat }_{\beta} \{(\beta  -\hat \beta )^2 =n\sigma^2/\sum_{i=1}^n x_i^2$, so its limit is a positive constant.
Thus, the model selector
\textit{achieves} the minimax rate
if and only if
\ben
\label{eq:minimaxrateattained}
\sup_\beta n \beta^2 P_{\beta}(A_n^c) \rightarrow \mbox {constant  as } n \rightarrow \infty.
\een
Yang (2005)  found a sequence  $\beta_n$ converging to zero, equivalently a sequence of  alternative models
converging to the null-model, \textit{slowly enough} to have $n \beta_n^2 \rightarrow \infty$,
 but at the same time \textit{fast enough} to ``confuse" the model selector,
   leading it  to choose the null model;  this keeps $P_{\beta_n}(A_n^c)$ away from zero, and actually arbitrarily close to one, so that  the minimax rate is not achieved.

\section{The Prediction/Minimax Rate Conflict}\label{sec:conflict}

In this section we describe the sequence $\beta_n$, chosen by Yang, to establish his result.  We also show that, with a minor modification, a similar sequence attains the minimax rate.  We then look a bit closer at Yang's sequence, showing how it  unfairly ``confuses" the model selector.

\subsection{Yang's Sequence}\label{sec:yang}
%
%
For the model selection problem described in (\ref{eq:simpleregressionmodel}), the model selector is \textit{consistent} if
$$
P(A_n) \rightarrow 0 \mbox{ if $H_0$ is true}\mbox{ and } P(A_n) \rightarrow1 \mbox{ if $H_1$ is true}.
$$
Yang (2005) shows that a consistent model selector cannot achieve the minimax prediction rate using the following argument.  He considers a sequence of UMP tests for the hypotheses (\ref{eq:simpleregressionmodel}) having   power function
\begin{equation}\label{eq:power}
\pi_n(\beta) =
\prob_{\beta} \{ \sommaxiyi \geq d_n \}
= \prob _{\beta} \{ Z \geq \frac{d_n-\beta \sommaxisq}{\sqrtsommaxisq} \},
\end{equation}
where $Z \sim N(0,1)$, and equates $\prob_{\beta=0}(A_n)$ with $\pi_n(0)$, for each $n$. Note that the set $A_n^\ast = \{(x_i,y_i): \sum_i x_iy_i > d_n\}$ defines a model selector based on the UMP test, where we would estimate $\beta_n$ with the least squares estimator if this set occurred.

The requirement of consistency means that we must have $\pi_n(0) \rightarrow 0$ which implies
\ben
\label{eq:dnsqrtsommaxisqRightarrowInfty}
\frac{d_n}{\sqrtsommaxisq} \rightarrow \infty.
\een
Yang's objective is to find a sequence of alternatives $\beta_n \rightarrow 0^+$
such that
\ben
\label{eq:C1C2}
C1:\quad n \beta_n^2 \rightarrow \infty \quad  \mbox{ and } \quad C2:  (1-\pi_n(\beta_n)) \rightarrow \mbox{constant}>0,
\een
which would imply that $n \beta_n^2 \prob_{\beta_n} (A_n^c)  \rightarrow \infty$ (because the test is UMP), and hence lead to the conclusion that a consistent model selector does \textit{not} attain the minimax rate, that is, it violates (\ref{eq:minimaxrateattained}).
Specifically, Yang's choice for $\beta_n$ is
\ben
\label{eq:yangbetan}
\beta_n=\frac{1}{2}\frac{d_n}{\sommaxisq},
\een
which satisfies $C1$ and $C2$.  (There is a  typo in Yang (2005), as noted by the author on his webpage: on lines 4 and 6 of page 947 the factor 2 should be in the denominator.)

We now look at  sequence (\ref{eq:yangbetan}) a bit more closely.  Write
\begin{equation}\label{eq:yangschoice}
\beta_n= \frac{c_n} {\sqrtsommaxisq}, \quad \mbox{with} \, c_n=\frac{1}{2}\frac{d_n}{\sqrtsommaxisq} \rightarrow \infty  \, ,
\end{equation}
and, recalling that $\sommaxisq \asymp n$,
we immediately conclude that we are  in scenario $1$ of Lemma \ref{lem:HellingerAffinity}, namely that the two hypotheses are strongly separated. In this case any type of test can occur, and with Yang's choice the power function is asymptotically zero:
\be
\pinbetan &=&\prob_{\beta_n} \{ \sommaxiyi \geq d_n \}
= \prob \{ Z \geq \frac{1}{2}\frac{d_n}{\sqrtsommaxisq}  \} \rightarrow 0,
\ee
which follows from  (\ref{eq:dnsqrtsommaxisqRightarrowInfty}) and the fat that  $\sum_i x_i y_i \sim \mbox{N}(\beta_n \sum_i x_i^2,\sum_i x_i^2)$ .

Yang's sequence for $\beta_n$ thus produces a worthless test, since both the size (by assumption) and the power (by construction) tend to zero.
In this case we know that there even exists a perfect sequence of tests (we will construct one below). Clearly Yang's argument holds for any sequence $\beta_n=b\frac{d_n}{\sommaxisq}$, $0<b < 1$. If $b=1$ we get
$\pinbetan \rightarrow 1/2$, which would still support Yang's argument, although the test is no longer worthless (but actually rather poor because its power is fixed at $1/2$ for each $n$).

Looking at Yang's result from a testing perspective reveals one of its weak points.  His result holds because he chooses a sequence of alternatives that, while producing  asymptotically a strong separation between the two models, converges to the null along a path leading to a worthless, or at best a mediocre, test. This happens despite the test being UMP, and despite the fact that, by a minor modification, one could get a perfect sequence of tests, as the following example shows.

\begin{example}
\label{ex:dnequal2times}
Consider the UMP test for the hypotheses (\ref{eq:simpleregressionmodel}). Choose
\ben
\label{eq:dnequal2times}
\beta_n=(1+b^{\prime})\frac{d_n}{\sommaxisq}, \, b^{\prime}>0.
\een
Then
\be
\pinbetan=\prob \{ Z \geq -b^{\prime}\frac{d_n}{\sqrtsommaxisq}  \} \rightarrow 1,
\ee
which follows from (\ref{eq:dnsqrtsommaxisqRightarrowInfty}) and the fact that $b^{\prime}>0$.  As a consequence, by choosing a sequence of alternatives which is structurally equivalent to Yang's,  although uniformly  larger by a factor $(1+b^{\prime})/b$,
$b^{\prime} >0$, $0<b<1$,
we get a perfect sequence of tests.

\cvd
\end{example}

With the choice of sequence (\ref{eq:dnequal2times}), $ 1-\pinbetan$ goes to zero exponentially fast,
and it is easy to show that $n \beta_n^2 ( 1-\pinbetan) \rightarrow 0$,  instead of going to $\infty$ as under Yang's choice (\ref{eq:yangbetan}). Thus, under this sequence of alternatives, we even beat the minimax rate!
Example \ref{ex:dnequal2times} reinforces the view that Yang's result  is based on a rather artificial sequence.
Next we provide further insight into his choice of sequence $\beta_n$.

\subsection{Confusing the Model Selector}

In Section \ref{sec:consistandpred} we noted that Yang's sequence $\beta_n$  was constructed in such a way as to ``confuse" the model selector.  We now look a bit more closely at this claim.

What happens with Yang's sequence is that, as $n \rightarrow \infty$, all of the mass of the distribution of $\sum_i x_i y_i$ is concentrated in the acceptance region of the test (that is, in $A_n^c$), even when $\beta_n >0$, making the ``correct" decision that of   accepting $H_0$.  To see this,
recall that, for $\beta=\beta_n$,  the UMP test-statistic $\sum_i x_i y_i$ is distributed as $\mbox{N}(\beta_n \sum_i x_i^2,\sum_i x_i^2)$. 
Consider the probability
$$
{\rm Pr}\left( \sum_i x_i y_i< \beta_n   \sum_i x_i^2 + M \sqrt{ \sum_i x_i^2}   \right),
$$
which grows arbitrarily close to 1 as $M$ increases. Now set $\beta_n=b \frac{d_n}{\sommaxisq}$, $0<b<1$ as in Yang's choice.
Then, for any $M>0$,
$
\beta_n  \sommaxisq + M \sqrtsommaxisq < d_n
$
eventually (that is, as $n$ grows), because the previous inequality  is equivalent to
\ben
\label{eq:M<(1-b)}
M < (1-b) \frac{d_n}{\sqrtsommaxisq},
\een
which holds true since  the right-hand-side tends to infinity because of (\ref{eq:dnsqrtsommaxisqRightarrowInfty}).
But this means that the support of the UMP-test statistic under this sequence of alternatives is eventually disjoint from the $H_1$-acceptance region
of the test: this is why the test is fooled and chooses $H_0$ incorrectly, with probability tending to one. Notice that if $b=1$ then (\ref{eq:M<(1-b)}) does not hold.

Finally, for the sequence with $\beta_n=(1+b^{\prime}) \frac{d_n}{\sommaxisq}$, $b^{\prime}>0$, then $
\beta_n  \sommaxisq - M \sqrtsommaxisq > d_n,$
eventually, because the previous inequality  is equivalent to
\ben
\label{eq:M>bprime}
-M > -b^{\prime} \frac{d_n}{\sqrtsommaxisq},
\een
which holds true since  the right-hand-side tends to $- \infty$ because of (\ref{eq:dnsqrtsommaxisqRightarrowInfty}).
This means that the support of the UMP-test statistic under this sequence of alternatives is eventually contained in the $H_1$-acceptance region of the test:
the test correctly chooses $H_1$ with probability tending to one, and thus chooses $H_0$ with probability tending to zero. The model-selector estimator attains the minimax rate.


\section{Contiguous Sequences}
\label{sec:contiguity}

\citet [Appendix C, Proposition C.1]{Leeb:Pots:2005}  present a result which  is comparable to that of Theorem 1 in \cite{Yang:2005}.
They consider  a linear regression model with mean structure $\alpha x_{1i}+\beta x_{2i}$ under the unrestricted case, and mean structure   $\alpha x_{1i}$ under the reduced model. They deal, among other things, with the scaled risk, under  squared-error loss, of the post-model selection  estimator (least squares) of $\alpha$.  As in \cite{Yang:2005}, they claim that its supremum diverges to infinity whenever the model selection procedure is consistent.
%
 Although Yang is concerned with prediction and not estimation, the connection between the two results is  apparent in the case of normal errors.
%
  Yet, Leeb and P\"{o}tscher's argument is quite different from Yang's,  because it relies on the notion of \textit{contiguity}. For an introduction to the notion of contiguity,  see \citet[sec. 6.2]{Vaar:1998} and  \citet[sec. 12.3]{Lehm:Roma:2005}.  Here we revisit Yang's problem using the notion of contiguity. We provide an evaluation of this technique for the problem at hand and raise some critical issues.

Let $P_n$ and $Q_n$ be measures on a measurable spaces $(\Omega_n, \mathcal{A}_n)$.
%
\begin{definition}
The sequence $Q_n$ is \textit{contiguous} with respect to the sequence $P_n$ if $P_n(A_n) \rightarrow 0$ implies $Q_n(A_n) \rightarrow 0$ for every sequence of measurable
sets $A_n$.
This is denoted $Q_n \lhd P_n$.
\end{definition}
One can regard contiguity as the asymptotic analogue of the classic notion of absolute continuity of measures. The strength of contiguity stems from the fact that $Q_n$-limit law of random vectors $U_n: \Omega_n \mapsto \Re^k$ can be obtained from suitable $P_n$-limit laws; the usefulness of such result is apparent when the latter calculations are much  easier than the former.
If $Q_n$ is contiguous with respect to  $P_n$, and \textit{viceversa}, then we write $Q_n \lhd \rhd P_n$.


The following result considers the model selection problem discussed by Yang, and relates it to the notion of contiguity.
%
\begin{prop}
\label{prop:Pnrsqrtncontiguous}
Consider the problem described in (\ref{eq:simpleregressionmodel}).
Let $\Pno$ be the sequence of probability measures under the null model $H_0$,
and $\Pnbetan$  be the  sequence of probability measures corresponding to the local alternative  models $H_{1n}
: y_i= \beta_n x_i + \epsilon_i, \, \beta_n>0$, $\beta_n \rightarrow 0^+$.
Then $\Pnbetan \lhd \Pno$ if and only if $\beta_n=O(1/\sqrt{n})$; additionally $\Pnbetan \lhd \rhd \Pno$ under the same condition.
\end{prop}

\begin{proof}
The proof is essentially the same as that for proving contiguity of the joint distribution of $n$ iid observations from a Normal  with mean $\xi_n$ and variance 1 with respect to the joint distribution of $n$ iid with observations  from a Normal with mean 0 and variance 1; see   \citet[examples 12.3.3 and  12.3.6]{Lehm:Roma:2005}.
To see why, simply notice  that the likelihood ratio is
$$d\Pnbetan/d\Pno=\exp \{ \beta_n \sommaxiyi-(\beta^2_n /2)\sommaxisq \}.$$
Under $\Pno$, $\sommaxiyi \sim N(0 , \sommaxisq)$,
and thus,
\be
\beta_n \sommaxiyi -(\beta^2_n/2) \sommaxisq
\sim N((-\beta^2_n/2) \sommaxisq,\beta^2_n \sommaxisq),
\ee
again under $\Pno$.
So the only difference  between the simple regression case we are discussing and the iid case from a Normal is that the former has $\sommaxisq$ while
the latter has $n$. Since these two quantities are asymptotically of the same order, one can use  the same  argument in either case.
\end{proof}

  From Proposition \ref{prop:Pnrsqrtncontiguous} it appears that $\Pnbetan \lhd \Pno$ if and only if the sequence
  $n \beta_n^2$ remains bounded; under the same condition mutual contiguity holds. In particular even if $\beta_n \rightarrow 0$, but at a slower rate than $1/\sqrt{n}$, as in Yang's case, see (\ref{eq:yangschoice}),
  then contiguity of $\Pnbetan$ fails.  Notice that contiguity  holds also  if $\beta_n=o(1/\sqrt{n})$, because $\beta_n^2 \sommaxisq$ goes to zero and hence is bounded. However this case is of no interest for proving failure to attain the minimax rate, because condition C1 in (\ref{eq:C1C2})  is not satisfied (that is, $n\beta_n^2$ does not diverge but actually goes to zero).

How can we use contiguity to obtain a result similar to Yang's? Here is the idea.
Yang's result obtains if we show that
\be
\limn n \beta_n^2 \prob_{\Pnbetan} \{ A^c_n \} = \infty
\ee
To exploit contiguity,
$\beta_n$ must be of order $1/\sqrt{n}$. Set for definiteness $\beta_n=r/\sqrt{n}$, for some  positive \textit{fixed} $r$.
We get
\be
\lim_{n \rightarrow \infty} \prob_{\Pnrsqrtn} \{ A^c_n \}=\lim_{n \rightarrow \infty} \prob_{\Pno} \{ A^c_n \}=1,
\ee
where the first equality sign follows from contiguity of $\Pnrsqrtn$  with respect  to $\Pno$, while the  second is a consequence of the assumed consistency of the model selector (recall that $A^c_n$ means
accepting $H_0: \beta=0$).
Therefore
\ben
\label{eq:limbetan2}
\limn n r^2(1 / \sqrt{n})^2 \prob_{\Pnrsqrtn} \{ A^c_n \}=r^2.
\een
At this stage it would seem that the minimax rate \textit{is} attained under this sequence, because the limit, however large, is finite.
To circumvent (\ref{eq:limbetan2}), and get the opposite conclusion that the 
rate is actually infinite,
one ought to apply  the argument in  \citet [p. 59]{Leeb:Pots:2005}, and let $r$ grow \textit{arbitrarily large} (technically this amounts  to take a further limit $r \rightarrow \infty$).
%
However, there is a subtle difference between this argument and Yang's result.

First of all, to conclude that $\limn n \beta_n^2 \prob_{\Pnbetan} \{ A^c_n \}=\infty$ one should prove that, for any $R>0$, there exists an $n_0$ such that
\ben
\label{eq:n>n0}
n >n_0 \Rightarrow n \beta_n^2 \prob_{\Pnbetan} \{ A^c_n \} >R.
\een
Having set $\beta_n=r/\sqrt{n}$, condition (\ref{eq:n>n0}) translates to
\ben
\label{eq:n>n0bis}
n>n_0 \Rightarrow \prob_{\Pnrsqrtn} \{ A^c_n \}>R/r^2.
\een
Since $r$ is \textit{fixed}, condition (\ref{eq:n>n0bis}) can be easily violated choosing for instance $R>r^2$.

Secondly, and possibly more importantly, the argument based on contiguity conveys the misconception that sequences of alternatives $\beta_n$ of order   $1/\sqrt{n}$
   can fail to attain the minimax rate. This would be quite surprising because, on the contrary, it is  well known  that this type of sequences is the \textit{only} one which makes sense for asymptotic comparison of tests, as lucidly remarked, for instance,  in \citet[example 12.3.6]{Lehm:Roma:2005}.

In the light of the above remarks, and of (\ref{eq:limbetan2}), it should be clear that sequences of alternatives  of order $1/\sqrt{n}$ \textit{do} achieve the minimax rate.

\section{Concluding Remarks}
\label{sec:concluding}
In this paper we have cast the  so-called AIC-BIC dilemma into perspective. On the one hand it is true that estimators and predictors based on consistent  model selection procedures may lead to an infinite scaled risk, thus failing to attain the usual minimax rate, as \cite{Yang:2005} showed.
On the other hand, this phenomenon occurs only for sequences of alternatives which are strongly separated from the null model
(this inflates the bias when the null model is chosen, while the alternative holds).
But such sequences are well known to be  of no use in asymptotic comparison of testing procedures, because they always admit a perfect sequence of tests (the power goes to 1 while the size goes to zero).

Additionally, the non-attainment of the minimax rate takes place only for a specific subset of  these sequences, namely those whose support under the alternative is eventually fully contained in the null-acceptance region. This explains why the problem occurs: the selector (not surprisingly!) chooses the null model, although the alternative holds.
Finally we have argued that contiguity arguments have little to say with regard to the AIC-BIC dilemma: contiguity is synonymous with sequences of alternatives converging to the null at the appropriate rate $1/\sqrt{n}$: no pathological behaviour can occur in this case.

\begin{center}

ACKNOWLEDGEMENT

\end{center}

 George Casella was supported by National Science Foundation Grants
DMS-04-05543, DMS-0631632 and SES-0631588.
Guido Consonni was  supported by MIUR PRIN 2007XECZ7L$_{-}$001.
This paper was begun while the second Author was visiting
the Department of Statistics, University of Florida,
Gainesville. Support and warm hospitality from this institution  is gratefully acknowledged.

\bibliographystyle{biometrika}
\bibliography{AICBIC}

\end{document}